\documentclass[12pt]{article}
\usepackage{latexsym}
\usepackage{hyperref}
\usepackage[pdftex]{graphicx}
\newcommand{\nc}{\newcommand}
\nc{\nt}{\newtheorem} \nt{thm}{Theorem}[section]
\nt{cor}[thm]{Corollary} \nt{prop}[thm]{Proposition}
\nt{lem}[thm]{Lemma} \nt{defn}[thm]{Definition}
\nt{exa}[thm]{Example} \nt{ass}[thm]{Assumption}
\nt{alg}[thm]{Algorithm} \nt{conj}[thm]{Conjecture}
\nt{com}[thm]{Comment for Dennis} \nt{coml}[thm]{Comment for Adrian}
\nc{\ip}[2]{\mbox{$\langle #1,#2 \rangle$}} \nc{\pf}{\noindent{\bf
Proof\ \ }} \nc{\finpf}{\hfill{$\Box$}\linespace}
\nc{\linespace}{\vspace{\baselineskip} \noindent} \nc{\bx}{\bar x}
\nc{\R}{{\bf R}} \nc{\Rn}{{\bf R}^n} \nc{\Sn}{{\bf S}^n}
\nc{\Rm}{{\bf R}^m} \nc{\E}{{\bf E}} \nc{\e}{\epsilon}
\nc{\rt}{\rightarrow} \nc{\Diag}{{\rm Diag}\,} \nc{\tra}{\mbox{\rm
tr}\,} \nc{\HH}{\mathcal{H}}
\def\tto{\;{\lower 1pt \hbox{$\rightarrow$}}\kern -12pt
           \hbox{\raise 2.8pt \hbox{$\rightarrow$}}\;}

\newenvironment{myequation}{\setcounter{equation}{\value{thm}}
   \begin{equation}}{\addtocounter{thm}{1}\end{equation}}

\nc{\bmye}{\begin{myequation}} \nc{\emye}{\end{myequation}}
\newcounter{algorithm}
\renewcommand{\thealgorithm}{\arabic{section}.\arabic{algorithm}}
\nc{\newalgorithm}{\stepcounter{algorithm}\thealgorithm}

\begin{document}

\title{Metric Subregularity and the Proximal Point Method}
\author{
D. Leventhal\thanks{School of Operations Research and Information
Engineering, Cornell University, Ithaca, NY 14853, U.S.A.
\texttt{leventhal\char64 orie.cornell.edu }}}

\maketitle

 \\
{\bf AMS 2000 Subject Classification: }

\begin{abstract}
We examine the linear convergence rates of variants of the proximal
point method for finding zeros of maximal monotone operators. We
begin by showing how metric subregularity is sufficient for linear
convergence to a zero of a maximal monotone operator. This result is
then generalized to obtain convergence rates for the problem of
finding a common zero of multiple monotone operators by considering
randomized and averaged proximal methods.
\end{abstract}
\section{Introduction}
Let $\mathcal{H}$ be a real Hilbert space and let
$T:\mathcal{H}\tto\mathcal{H}$ be a set-valued mapping. Two common
problems that arise in several branches of applied mathematics are
to \bmye\label{Problem1}\mbox{Find }x\in\mathcal{H}\mbox{~~~~such
that~~}0\in T(x)\emye and, more generally
\bmye\label{Problem2}\mbox{Find }x\in\mathcal{H}\mbox{~~~~such
that~~}0\in\cap_{i\in I}T_i(x),\emye where $I$ is some index set.
Specifically, these problems correspond to finding a zero of an
operator and, more generally, a common zero of multiple operators.

Suppose that the operators under consideration are {\em monotone},
meaning that $$\langle x_1-x_0, y_1-y_0\rangle \geq 0 \mbox{ for all
}x_0,x_1\in \mathcal{H}, y_0\in T(x_0), y_1\in T(x_1).$$ For
$\lambda
> 0$, the mappings $J_{\lambda T}:= (I+\lambda T)^{-1}$ are the {\em
resolvents} of $T$, which were shown to be at most single-valued in
\cite{Minty62}. One proposed method for solving Problem
\ref{Problem1} is the {\em proximal point algorithm}, considered
originally in \cite{Martinet70} and more thoroughly explored by
\cite{Rockafellar76}, given by, for $k=0,1,2,\ldots,$
\bmye\label{proxpoint}x_{k+1} = J_{\lambda T}(x_k).\emye

Our goal is to examine how appropriate regularity assumptions on the
operators $T$ (or $T_1,\ldots,T_m$, respectively) affect the speed
of convergence of variants of the proximal point algorithm. In order
to do so, the remainder of this paper is organized as follows. In
Section \ref{background}, we provide notation and basic facts about
monotone operators, metric regularity and subregularity, and the
geometry of convex sets. Then, in Section \ref{convergence}, we show
how assumptions of metric subregularity can be used to demonstrate
linear convergence of both the proximal point algorithm for Problem
\ref{Problem1} and a randomized proximal point algorithm for Problem
\ref{Problem2}.

\section{Background and Notation}\label{background}
A single-valued operator $U$ is {\em firmly non-expansive} if
\bmye\label{fne}\|U(x)-U(y)\|^2 + \|(I-U)(x) - (I-U)(y)\|^2 \leq
\|x-y\|^2 ~~\forall x,y\in \mathcal{H}\emye It was shown in
\cite{Rockafellar76, Eckstein89} that an operator $T$ is monotone if
and only if its resolvents are firmly non-expansive. The domain of
$T$ is $\{x\in\mathcal{H}: T(x)\neq\emptyset\}$ and the inverse
operator, $T^{-1}$, is defined by $T^{-1}(y) = \{x: y\in T(x)\}$. It
is known that (see \cite{Rockafellar98}, for example) $T$ is
monotone if and only if $T^{-1}$ is monotone and, if $T$ is {\em
maximal monotone}, meaning the graph of $T$ is not strictly
contained in the graph of another monotone operator, then both $T$
and $T^{-1}$ are closed and convex-valued and the domain of the
resolvents of $T$ is $\mathcal{H}$.

We are interested in how certain regularity conditions affect {\em
local} rates of convergence. One prominent condition is the idea of
metric regularity of set-valued mappings. We say the set-valued
mapping $\Phi$ is \textit{metrically regular} at $\bar x$ for $\bar
b\in \Phi(\bar x)$ if there exists $\gamma > 0$ such that \bmye
\label{MetReg} d(x, \Phi^{-1}(b)) \leq \gamma~ d(b, \Phi(x)) \mbox{
for all } (x,b) \mbox{ near } (\bar x, \bar b). \emye Further, the
\textit{modulus of regularity} is the infimum of all constants
$\gamma$ such that Inequality \ref{MetReg} holds.

A slightly weaker condition is that of metric subregularity. We say
the set-valued mapping $\Phi$ is \textit{metrically subregular} at
$\bar x$ for $\bar b\in \Phi(\bar x)$ if there exists $\gamma > 0$
such that \bmye \label{MetSubreg} d(x, \Phi^{-1}(\bar{b})) \leq
\gamma~ d(\bar{b}, \Phi(x)) \mbox{ for all } x \mbox{ near } \bar x.
\emye Further, the \textit{modulus of subregularity} is the infimum
of all constants $\gamma$ such that Equation \ref{MetSubreg} holds.
Note that for metric subregularity, the reference vector $\bar{b}$
is fixed in Inequality \ref{MetSubreg} but not in Inequality
\ref{MetReg}. It is clear from the definitions that metric
regularity implies metric subregularity; hence, the modulus of
subregularity is no larger than the modulus of regularity, using the
convention that the modulus of (sub)regularity is infinite if the
mapping fails to be metrically (sub)regular.

The property of metric regularity is connected with other ideas in
variational analysis. The simplest connection, as shown in \cite[Ex.
1.1]{Dontchev}, is that metric regularity generalizes the Banach
open mapping principle, essentially saying that a bounded and linear
mapping is metrically regular if and only if it is surjective; in
such a case, the modulus of regularity is simply $\sup_{y\in
B}\{d(0,A^{-1}(y))\}$ where $B$ is the unit ball. If the mapping
$\Phi$ has a closed-convex graph, the Robinson-Ursescu Theorem says
that $\Phi$ is metrically regular at $\bar{x}$ for $\bar{y}$ if and
only if $\bar{y}$ is in the interior of the range of $\Phi$. Metric
regularity is also known to be equivalent to several others in
variational analysis, namely the Aubin property of $\Phi^{-1}$ and
the openness at linear rate of $\Phi$. Additionally, metric
regularity has been shown to be a generalization of the Eckart-Young
result from matrix analysis on the distance to singularity of a
matrix. Further, a result originating with Lyusternik and Graves
(\cite{Lyusternik34}, \cite{Graves50}) and extended by others (for
example, \cite{Dmitruk80},\cite{Ioffe00}, \cite{Dontchev}) show that
metric regularity is determined by the first-order behavior of a
mapping and is preserved by sufficiently small first-order
perturbations. Additional information about metric regularity and
its relationship to other concepts in variational analysis can be
found in \cite{Dontchev04}, \cite{Ioffe00}, and \cite{Dontchev},
among others.

A central tool frequently appearing in variational analysis is that
of the {\em normal cone} of a closed, convex set, $S$. Specifically,
the normal cone of $S$ at $\bar{x}\in S$ can be defined as
\bmye\label{normalcone}N_{S}(\bar{x}) := \{x^*\in
\mathcal{H}:~\langle x^*, s-\bar{x}\rangle \leq 0~~\forall s\in
S\}\emye and $N_S(\bar{x}) = \emptyset$ if $\bar{x}\not\in S$. Let
$d(x,S)$ denote the distance from $x$ to $S$, given by $d(x,S) :=
\inf_{s\in S}\|x-s\|$. Further, let $P_S(x)$ be the projection
operator onto $S$, i.e., the set of such minimizers. If $S$ is
closed, convex and non-empty, then $P_S$ is single-valued
everywhere. Further, the projection operator is firmly non-expansive
(\cite[Thm 5.5]{Deutsch}) and can be characterized by
\bmye\label{projchar}z = P_S(x) \leftrightarrow z\in S \mbox{~and~}
x-z\in N_S(z).\emye

A method of characterizing regularity of closed sets
$S_1,\ldots,S_m$ is by considering regularity properties of a
related set-valued mapping. Given a Hilbert space, $\HH$, consider
the product space $\HH^m$ with the induced inner product defined by
$$\langle (x_1,x_2,\ldots,x_m), (y_1,y_2,\ldots,y_m)\rangle = \sum_{i=1}^m \langle x_i,y_i\rangle$$
and consider the set-valued mapping given by $\Phi(x) =
[S_1-x,\ldots,S_m-x]^T$. Note that $0\in\Phi(x)$ if and only if
$x\in\cap_i S_i$. Using metric regularity as a starting point,
suppose $\Phi(x)$ is metrically regular at $\bar{x}$ for 0. From the
definition, metric regularity of $\Phi$ at $\bar{x}$ for 0 is
equivalent to the {\em strong metric inequality}, examined in
\cite{Kruger05} and \cite{Kruger06}, among others, defined by the
existence of $\beta,\delta>0$ such that, for $i=1,\ldots,m$,
\bmye\label{strongmetric}\mbox{~~~~~}d(x, \cap_i (S_i - z_i)) \leq
\beta \max_{1\leq i\leq m}d(x+z_i, S_i)\mbox{~~for all~}x\in
\bar{x}+\delta B,\mbox{~}z_i\in\delta B.\emye Characterizing this in
terms of normal cones, it was shown in \cite[Thm. 1, Prop. 10, Cor.
2]{Kruger06} that this is equivalent to the existence of a constant
$k>0$ such that \bmye\label{quantreg}~~z_i\in\delta B,~y_i\in
N_{S_i}(\bar{x}+z_i)~~(i=1,\ldots,m) \Rightarrow \sum_i\|y_i\|^2
\leq k^2\|\sum_iy_i\|^2.\emye By using the formula in \cite[Thm
9.43]{Rockafellar98} for expressing the modulus of regularity in
terms of coderivatives, it was shown in \cite{Lewis08} that the
modulus of regularity of $\Phi$ at $\bar{x}$ for 0 equals
$$\lim_{\delta\downarrow0}\Big\{\inf\{k:
\mbox{~Inequality~\ref{quantreg} holds.}\}\Big\},$$ with this value
being infinite being equivalent to a lack of metric regularity of
$\Phi$.

Consider a relaxed variant of the strong metric inequality, known
simply as the {\em metric inequality} as studied in \cite{Ioffe00},
\cite{Ngai01} and \cite{Kruger06} among others, defined to hold at
$\bar{x}$ if there exists $\beta>0$ such that
\bmye\label{metricineq}d(x,\cap_i S_i) \leq \beta \max_{1\leq i\leq
m} d(x, S_i)\mbox{~~for all~}x\in \bar{x}+\delta B.\emye If
Inequality \ref{metricineq} is valid for $\delta = \infty$, we
obtain the property of linear regularity and if it holds for all
$\delta > 0$, it is equivalent to the property of bounded linear
regularity, as studied in \cite{Bauschke}, \cite{Bauschke96},
\cite{Bauschke97}, \cite{Bauschke00}, \cite{Beck03} and others,
often in an algorithmic context. It is easy to show that the
existence of a $\delta>0$ such that Inequality \ref{metricineq}
holds is equivalent to the previously defined mapping $\Phi$ being
metrically subregular at $\bar{x}$ for 0.

Our focus for the remainder of this paper will involve metric
subregularity. Unfortunately, several of the stability properties
and some of the geometric intuition that accompanies metric
regularity---especially that relating to normal cones of
sets---fails to have a natural equivalent for metric subregularity;
some examples of this phenomenon are given in \cite{Dontchev04}.
However, since metric regularity implies metric subregularity, the
intuition provided by metric regularity can be applied to the
following results when that property does, in fact, hold.
Additionally, if the monotone operators under consideration are
actually subdifferentials of convex functions, characterization of
both metric regularity and subregularity  in terms of the underlying
function was shown in \cite{Artacho08}, providing additional
intuition.

\section{Metric Regularity and Linear Convergence}\label{convergence}
We now return to Problem \ref{Problem1}, the problem of finding a
zero of a maximal monotone operator. Variants of proximal point
algorithms for solving this and related problems have been
considered by a wide variety of authors, including
\cite{Rockafellar76},
\cite{Luque84},\cite{Solodov99},\cite{Pennanen02}, \cite{Artacho07}
and others.

Many authors consider an algorithmic framework much more general
than the one considered in this paper. Some of the better-studied
variants allow for a varying proximal parameter $\lambda$, allow
approximate computation of the proximal iteration, allow over- or
under-relaxation in the proximal step or incorporate an additional
projective framework. These ideas have often proven worthwhile both
for designing a computationally practical and efficient algorithm as
well as for improving the convergence analysis. However, in this
paper, we will only consider algorithms in their ``classical'' form,
assuming exact computation of the resolvent with a fixed proximal
parameter. Our particular interest is in exploring how naturally
occurring constants---for example, the modulus of subregularity of
the mappings themselves and of the mapping associated with the
solution sets---govern the local rate of convergence and, further,
how randomization as an analytical tool can emphasize this
connection. To begin, consider the basic proximal point algorithm
given by \ref{proxpoint}, where $x_{k+1} = J_{\lambda T}(x_k)$.
Under an assumption of metric subregularity, we obtain the following
initial result.

\begin{thm}\label{lconvprox1}Suppose $T$ is maximal monotone and metrically subregular at
$\bar{x}\in T^{-1}(0)$ for 0 with regularity modulus $\gamma$. Let
$\bar{\gamma}
> \gamma$ and suppose $x_0$ is sufficiently near $\bar{x}$. Then the iterates given by Algorithm \ref{proxpoint}
are linearly convergent to $T^{-1}(0),$ the zero-set of $T$,
satisfying
$$d(x_{k+1}, T^{-1}(0))^2 \leq \frac{\bar{\gamma}^2}{\lambda^2 +
\bar{\gamma}^2}d(x_k,T^{-1}(0))^2.$$ \end{thm}

\pf Let $\hat{x}\in T^{-1}(0)$ and note that $J_{\lambda T}(\hat{x})
= \hat{x}.$ Since the resolvent of a monotone operator is firmly
non-expansive, it follows from Inequality \ref{fne} that, for any
$x$,
$$\|J_{\lambda T}(x) - J_{\lambda T}(\hat{x})\|^2 \leq
\|x-\hat{x}\|^2 - \|(I-J_{\lambda T})(x) - (I-J_{\lambda
T})(\hat{x})\|^2,$$ implying that \bmye\label{lcproof1}\|J_{\lambda
T}(x)-\hat{x}\|^2 \leq \|x - \hat{x}\|^2 - \|x - J_{\lambda
T}(x)\|^2.\emye However, by definition of $J_{\lambda T}$,
$$x- J_{\lambda T}(x) \in \lambda T(J_{\lambda T}(x)).$$ In
particular,

\bmye\label{stepbound}\|x- J_{\lambda T}(x)\|  \geq  \lambda
\min\{\|z\|: z\in T(J_{\lambda T}(x))\}  = \lambda~d(0,T(J_{\lambda
T}(x))).\emye Now, note that since the resolvents and projection
operators are firmly non-expansive, if $x_0$ has the property of
being sufficiently close to $\bar{x}$ such that Inequality
\ref{MetSubreg} holds with constant $\bar{\gamma}$, then $x_j$ and
$P_{T^{-1}(0)}(x_j)$ do as well for each $j\geq 0$. Therefore, it
follows that
\begin{eqnarray*}
\lefteqn{d(x_{k+1}, T^{-1}(0))^2}&\hspace{25pt}&\\
 & \leq &
\|x_{k+1}-P_{T^{-1}(0)}(x_k)\|^2\hspace{25pt}\\
\hspace{-15pt}& \leq & \|x_k - P_{T^{-1}(0)}(x_k)\|^2 - \|x_k - J_{\lambda T}(x_k)\|^2 \mbox{~~~(Inequality~\ref{lcproof1})}\\
& \leq & d(x_k, T^{-1}(0))^2 - \lambda^2d(0, T(J_{\lambda
T}(x_k)))^2 \mbox{~~~(Inequality \ref{stepbound})}\\
& \leq & d(x_k, T^{-1}(0))^2 -
\frac{\lambda^2}{\bar{\gamma}^2}d(J_{\lambda T}(x_k), T^{-1}(0))^2
\mbox{~~~(Inequality \ref{MetSubreg})}\\
& = &d(x_k, T^{-1}(0))^2 -
\frac{\lambda^2}{\bar{\gamma}^2}d(x_{k+1}, T^{-1}(0))^2.
\end{eqnarray*}
This implies that
$$(1+\frac{\lambda^2}{\bar{\gamma}^2})d(x_{k+1}, T^{-1}(0))^2 \leq
d(x_k, T^{-1}(0))^2,$$ from which the result follows.\finpf

Further observe that by considering a sequence $\{\lambda_k\}$ such
that $\lambda_k\rightarrow\infty$ instead of a fixed $\lambda$ in
the above algorithm, we obtain superlinear convergence.

Our primary interest in Theorem \ref{lconvprox1} is as a tool in
proving the following result, Theorem \ref{randomproxlc}. However,
we note that Theorem \ref{lconvprox1} is similar to some previously
known results. For example, linear convergence was shown in
\cite{Rockafellar76} and \cite{Solodov99}, under a framework that
permitted error in evaluating the resolvent, with a slightly
stronger regularity assumption. In particular, as a limiting case
with no such error in evaluating the resolvent, an identical
convergence rate was obtained in \cite{Rockafellar76}. The result by
Solodov and Svaiter in \cite{Solodov99}, however, corresponds to a
hybrid proximal-projection algorithm.

We wish to generalize this result to that of Problem \ref{Problem2},
finding a common zero among a group of maximal monotone operators,
$T_1,\ldots,T_m$. Variants of proximal point algorithms for this
problem have been considered by a variety of authors, including
\cite{Kiwiel97}, \cite{Lehdili99}, \cite{Combettes04},
\cite{Hirstoaga06}, among others. In what follows, consider the
following randomized variant of a proximal point algorithm: for
$k=0,1,2,\ldots,$ \bmye\label{proxpoint2}x_{k+1} = J_{\lambda
T_i}(x_k)\mbox{~~~~with probability~~}\frac{1}{m},
~~~i=1,\ldots,m.\emye Then we obtain the following result.
\begin{thm}\label{randomproxlc}Suppose the following assumptions hold:
\begin{enumerate}
\item The maximal monotone operators $T_i$, $i=1,\ldots,m$, are
metrically subregular at $\bar{x}\in\cap_jT_j^{-1}(0)$ for 0 with
modulus $\gamma_i$.
\item The mapping $\Phi(x) = [T_1^{-1}(0)-x, \ldots, T_m^{-1}(0)-x]^T$ is metrically
subregular at $\bar{x}$ for 0 with modulus $\kappa$.
\item $\bar{\gamma} > \max\{\gamma_1,\ldots,\gamma_m\}$ and $\bar{\kappa} > \kappa$.
\item $\lambda^2 > 3\bar{\gamma}^2$.
\end{enumerate}
Then for $x_0$ sufficiently close to $\bar{x}$, Algorithm
\ref{proxpoint2} satisfies
$$d(x_{k+1}, \cap_j T_j^{-1}(0)) \leq d(x_k, \cap_j
T_j^{-1}(0))$$ and
$$\E[d(x_{k+1}, \cap_j T_j^{-1}(0))^2 ~|~ x_k] \leq \Big(1-\frac{1}{m\bar{\kappa}^2} +
\frac{2}{m\bar{\kappa}^2}\Big(\frac{\bar{\gamma}^2}{\lambda^2+\bar{\gamma}^2}\Big)^{\frac{1}{2}}\Big)d(x_k,
\cap_j T_j^{-1}(0))^2.$$
\end{thm}

\pf If $x_0$ is sufficiently close to $\bar{x}$ such that Inequality
\ref{MetSubreg} holds with constant $\bar{\gamma}$ for each mapping
$T_i$, it follows from the firm non-expansivity of the resolvents
and the projection operator that each iterate $x_k$ and the
projection of each iterate onto the common zero set, $P_{\cap_j
T_j^{-1}(0)}(x_k)$, are sufficiently close to $\bar{x}$ as well.
Additionally, this implies the first conclusion of the theorem.

Suppose that at iteration $k$, the resolvent $J_{\lambda T_i}$ is
chosen by the algorithm. Then it follows that
\begin{eqnarray*}
\lefteqn{d(J_{\lambda T_i}(x_k), \cap_j T_j^{-1}(0))^2}
&\hspace{25pt}&\\ & = & \|J_{\lambda
T_i}(x_k) - P_{\cap_j T_j^{-1}(0)}(J_{\lambda T_i}(x_k))\|^2\\
& \leq & \|J_{\lambda T_i}(x_k) - P_{\cap_j
T_j^{-1}(0)}(x_k)\|^2 \hspace{40pt}\mbox{~~(Definition of Projection)}\\
& \leq & d(x_k,\cap_j T_j^{-1}(0))^2 -
\|x_k - J_{\lambda T_i}(x_k)\|^2 \hspace{10pt}\mbox{~~~~(Inequality \ref{lcproof1})} \\
& = & d(x_k,\cap_j T_j^{-1}(0))^2 - \Big\|\Big[x_k -
P_{T_i^{-1}(0)}(x_k)\Big] +
\Big[P_{T_i^{-1}(0)}(x_k) - J_{\lambda T_i}(x_k)\Big]\Big\|^2\\
& \leq & d(x_k,\cap_j T_j^{-1}(0))^2 - d(x_k, T_i^{-1}(0))^2 -
\|P_{T_i^{-1}(0)}(x_k) - J_{\lambda T_i}(x_k)\|^2 \\
& & ~~-~2\langle x_k - P_{T_i^{-1}(0)}(x_k), P_{T_i^{-1}(0)}(x_k) -
J_{\lambda T_i}(x_k)\rangle.
\end{eqnarray*}
Note that
\begin{eqnarray*}
\lefteqn{-2\langle x_k - P_{T_i^{-1}(0)}(x_k), P_{T_i^{-1}(0)}(x_k)
-
J_{\lambda T_i}(x_k)\rangle} \hspace{50pt}\\
& = & 2\langle x_k - P_{T_i^{-1}(0)}(x_k), [J_{\lambda T_i}(x_k) -
P_{T_i^{-1}(0)}(J_{\lambda T_i}(x_k))\rangle\\
&& \hspace{15pt}+\mbox{~~}\langle x_k - P_{T_i^{-1}(0)}(x_k),
P_{T_i^{-1}(0)}(J_{\lambda T_i}(x_k)) -
P_{T_i^{-1}(0)}(x_k)\rangle \\
& \leq & 2\langle x_k - P_{T_i^{-1}(0)}(x_k), J_{\lambda T_i}(x_k) -
P_{T_i^{-1}(0)}(J_{\lambda T_i}(x_k))\rangle \\
& \leq & 2\Big\|x_k - P_{T_i^{-1}(0)}(x_k)\Big\|~\Big\|J_{\lambda
T_i}(x_k) -
P_{T_i^{-1}(0)}(J_{\lambda T_i}(x_k))\Big\|\\
& = & 2~d(x_k, T_i^{-1}(0))~d(J_{\lambda T_i}(x_k),
T_i^{-1}(0))\\
& \leq &
2\Big(\frac{\bar{\gamma}^2}{\lambda^2+\bar{\gamma}^2}\Big)^{\frac{1}{2}}d(x_k,
T_i^{-1}(0))^2.
\end{eqnarray*}
The first inequality comes from the fact that
\linebreak$x_k-P_{T_i^{-1}(0)}(x_k)\in
N_{T_i^{-1}(0)}(P_{T_i^{-1}(0)}(x_k))$ so Inequality
\ref{normalcone} can be applied from the definition of the normal
cone. The second inequality is an application of the Cauchy-Schwartz
inequality. The rest follows from the definition of the projection
operator, followed by applying Theorem \ref{lconvprox1}, the
previous linear convergence result. Putting this together, we obtain
$$d(J_{\lambda T_i}(x_k), \cap_j T_j^{-1}(0))^2 \leq d(x_k, \cap_j T_j^{-1}(0))^2
-
\Big(1-2\Big(\frac{\bar{\gamma}^2}{\lambda^2+\bar{\gamma}^2}\Big)^{\frac{1}{2}}\Big)d(x_k,
T_i^{-1}(0))^2.$$ Taking the expected value, we obtain
\begin{eqnarray*}
\lefteqn{\E[d(x_{k+1}, \cap_j T_j^{-1}(0))^2 ~|~x_k]}\hspace{50pt}\\
& \leq & d(x_k, \cap_j T_j^{-1}(0))^2 -
\frac{1}{m}\Big(1-2\Big(\frac{\bar{\gamma}^2}{\lambda^2+\bar{\gamma}^2}\Big)^{\frac{1}{2}}\Big)\sum_{i=1}^md(x_k,
T_i^{-1}(0))^2\\
& = & d(x_k, \cap_j T_j^{-1}(0))^2 -
\frac{1}{m}\Big(1-2\Big(\frac{\bar{\gamma}^2}{\lambda^2+\bar{\gamma}^2}\Big)^{\frac{1}{2}}\Big)d(0,
\Phi(x_k))^2\\
& \leq & \Big(1-\frac{1}{m\bar{\kappa}^2} +
\frac{2}{m\bar{\kappa}^2}\Big(\frac{\bar{\gamma}^2}{\lambda^2+\bar{\gamma}^2}\Big)^{\frac{1}{2}}\Big)d(x_k,
\cap_j T_j^{-1}(0))^2,
\end{eqnarray*}
where the last inequality follows from the metric subregularity of
the mapping $\Phi(x) = [T_1^{-1}(0)-x, \ldots, T_m^{-1}(0) - x]^T$.
\finpf

Note that the last assumption in Theorem \ref{randomproxlc} is so
that the convergence rate is less than 1. Additionally, this type of
convergence result implies that $d(x_k, \cap_j
T_j^{-1}(0))\rightarrow 0$ almost surely (cf.\ \cite{Leventhal08}).

One particularly simple way of de-randomizing Algorithm
\ref{proxpoint2} is by considering averaged resolvents or, in the
terminology of \cite{Lehdili99}, the {\em barycentric proximal
method}. Specifically, given maximal monotone operators
$T_i,~i=1,\ldots,m$ with respective resolvents $J_{\lambda
T_i},~i=1,\ldots,m$, consider the algorithm described such that, for
$k=0,1,2,\ldots,$ \bmye\label{baryprox} x_{k+1} =
\frac{1}{m}\sum_{i=1}^m J_{\lambda T_i}(x_k) \emye and the
associated fixed-point problem \bmye\label{Problem3}\mbox{Find
}x\in\mathcal{H}\mbox{~~~~such that~~}x = \frac{1}{m}\sum_{i=1}^m
J_{\lambda T_i}(x).\emye The following proposition, found in
\cite{Lehdili99}, provides the necessary connection.
\begin{prop}[\cite{Lehdili99}]If $\bar{x}\in\cap_i T_i^{-1}(0)$, then $\bar{x}$ is a solution to Problem
\ref{Problem3}. Further, if $\cap_i T_i^{-1}(0)\neq \emptyset$, the
fixed points of Problem \ref{Problem3} are common fixed points of
all the $T_i$'s.
\end{prop}

Considering the example where each operator $T_i$ is the normal cone
mapping for some closed, convex set, it follows that Algorithm
\ref{baryprox} is simply the {\em averaged projections algorithm}
studied by \cite{Pierra}, \cite{Reich83}, \cite{Bauschke},
\cite{Lewis08}, and \cite{Leventhal08}, among others. More
generally, we can use the result of Theorem \ref{randomproxlc} to
generalize a result on averaged projections found in \cite[Thm
5.8]{Leventhal08} to the barycentric proximal method.

\begin{thm}Suppose the assumptions of Theorem \ref{randomproxlc}
hold. Then the conclusions of Theorem \ref{randomproxlc} hold for
Algorithm \ref{baryprox} as well.
\end{thm}
\pf Let $x_k$ be the current iterate, $x_{k+1}^{BP}$ be the new
iterate in the barycentric proximal method, Algorithm
\ref{baryprox}, and let $x_{k+1}^{RP}$ be the new iterate in the
randomized proximal point method, Algorithm \ref{proxpoint2}. First,
note that since each set $T_i^{-1}(0)$ is convex, the distance
function $d(~\cdot~, \cap_j T_j^{-1}(0))$ is as well, and
$$d(J_{\lambda T_i}(x_k),\cap_j T_j^{-1}(0)) \leq d(x_k, \cap_j
T_j^{-1}(0)) \mbox{~~for~} i=1,\ldots,m,$$ from which it follows
that
$$d(x_{k+1}^{BP}, \cap_j T_j^{-1}(0)) \leq d(x_k, \cap_j
T_j^{-1}(0)).$$ Let $\alpha = \Big(1-\frac{1}{m\bar{\kappa}^2} +
\frac{2}{m\bar{\kappa}^2}\Big(\frac{\bar{\gamma}^2}{\lambda^2+\bar{\gamma}^2}\Big)^{\frac{1}{2}}\Big)$
and observe that the function $d(~\cdot~,\cap_j T_j^{-1}(0))^2$ is
also convex. Noting that
\[
x_{k+1}^{BP} = \frac{1}{m}\sum_{j=1}^m J_{\lambda T_j}(x_k) =
\E[x_{k+1}^{RP} ~|~x_k],
\]
it follows that
\begin{eqnarray*}
d(x_{k+1}^{BP}, \cap_j T_j^{-1}(0))^2 & = & d(\E[x_{k+1}^{RP}
~|~x_k], \cap_j T_j^{-1}(0))^2\\ & \leq & \E[d(x_{k+1}^{RP}, \cap_j
T_j^{-1}(0))^2 ~|~ x_k]\\ & \leq & \alpha d(x_k, \cap_j
T_j^{-1}(0))^2,
\end{eqnarray*} from an application of Jensen's Inequality. \finpf

In particular, the barycentric proximal method converges at least as
quickly as the randomized proximal point method.

\bibliographystyle{plain}
\bibliography{dennis-references}

\end{document}